\documentclass[a4paper]{amsart}

\usepackage{color}
\usepackage{graphicx}
\usepackage{verbatim}
\usepackage[latin1]{inputenc}
\usepackage[T1]{fontenc}
\usepackage{amsthm}
\usepackage{amsmath}
\usepackage{amssymb}
\usepackage{relsize}
\usepackage[colorlinks=true,linkcolor=blue,urlcolor=blue,citecolor=blue]{hyperref}

\usepackage [all] {xy}
\usepackage{mathrsfs}

\newtheorem{theorem}{Theorem}[section]

\newtheorem{proposition}[theorem]{Proposition}
\newtheorem*{theorem*}{Theorem 1.13}

\theoremstyle{definition}
\newtheorem{definition}[theorem]{Definition}
\newtheorem{remark}[theorem]{Remark}
\newtheorem{example}[theorem]{Example}

\newtheorem{notation}[theorem]{Notation}
\newtheorem{question}[theorem]{Question}

\newcommand{\het}{H^2_{\text{\'et}}(X_{\overline K}, \fQ_l)}
\newcommand\blfootnote[1]{%
		\begingroup
		\renewcommand\thefootnote{}\footnote{#1}%
		\addtocounter{footnote}{-1}%
		\endgroup
		}

\newcommand{\sch}{\textbf{Sch}}

\newcommand{\spec}{\text{Spec}}

\newcommand{\id}{\text{Id}}

\newcommand{\gal}{\text{Gal}}

\newcommand{\fA}{\mathbb A}
\newcommand{\fC}{\mathbb C}

\newcommand{\fL}{\mathbb L}
\newcommand{\fP}{\mathbb P}
\newcommand{\fQ}{\mathbb Q}

\newcommand{\fZ}{\mathbb Z}

\newcommand{\sP}{\mathcal P}

\newcommand{\sO}{\mathcal{O}}

\newcommand{\sA}{\mathcal A}

\newcommand{\sM}{\mathcal M}

\newcommand{\sX}{\mathcal X}
\newcommand{\sY}{\mathcal Y}

\author{Author: Luigi Lunardon}
\title{What is the monodromy property for degenerations of Calabi-Yau varieties?}
\date{\today}
\address{Luigi Lunardon\\Office 671\\ Huxley Building\\ 180 Queen's Gate \\ SW7 2AZ\\Kensington\\ London.}
\email{l.lunardon16@imperial.ac.uk}

\begin{document}
\begin{abstract}
In this survey we discuss the state of art about the monodromy property for Calabi-Yau varieties. We explain what is the monodromy property for Calabi-Yau varieties and then discuss some examples of Calabi-Yau varieties that satisfy this property. After this recap, we discuss a possible approach to future research in this area.
\end{abstract}
\maketitle
\blfootnote{\textup{2010} \textit{Mathematics Subject Classification}: Primary: 14E18; Secondary: 14G10, 14D06.}
\section*{Introduction}
Our aim is the study of the monodromy of degenerations of Calabi-Yau varieties, and we are interested, in particular, in the so-called motivic monodromy conjecture. While degenerations and monodromy are intuitive concepts in complex geometry, translating these ideas in the setting of algebraic geometry is not trivial. In this section we try to give an intuitive description of our topic; hopefully, this picture will guide the reader through the understanding of this paper. 

Denote by \(D\subset \fC\) the set \(\{x\in \fC\colon |x|<1\}\). For the moment, a model for a degeneration of a Calabi-Yau is the datum of: a Calabi-Yau variety \(X\) over \(\fC\), a smooth complex variety \(\sX\), endowed with a proper fibration \(\pi\colon \sX\to D\) which is smooth on \(D^\ast=D\setminus\{0\}\) and satisfying the following additional properties:
\begin{enumerate} 
\item there exists \(z_0\in D, \; z_0\neq 0\) such that  \(\pi^{-1}(z_0) = X_{z_0} \cong X\);
\item \(X_0\) is a strict normal crossing divisor on \(\sX\) - by this we mean that the irreducible components of \((X_0)_{red}\) are smooth and with transverse intersection.
\end{enumerate}
The divisor \(X_0 = \sum_{i\in I} N_i E_i\) contains some of the geometrical data of the degeneration, but it is not enough to determine the motivic zeta function. As we see in Section \ref{2} some additional information is required.

Up to homotopy, we have a natural action of \(\fZ = \pi_1(D^\ast)\) on the underlying topological space of \(X\) This action induces an action on the cohomology groups \(H^i(X, \fQ)\). We can think of this as the cohomological datum of the degeneration. 

A natural question is how the geometrical and the cohomological data of a degeneration are related. The motivic monodromy conjecture suggests a possible answer to this question; this conjecture states that poles of the motivic zeta function of a degeneration of Calabi-Yau varieties are linked to monodromy eigenvalues. 
The motivic zeta function is a power series which depends on \(\sX^\ast\) - the restriction of the model \(\sX\) over \(D^\ast\) - and on a relative volume form \(\omega\) on \(X^\ast\). 
The motivic zeta function encodes all these data and some more information about the geometry of the special fiber. A degeneration of \(X\) satisfies the monodromy property if any pole of the motivic zeta function determines an eigenvalue in monodromy. The motivic monodromy conjecture states that Calabi-Yau varieties satisfy the monodromy property.

The motivic monodromy conjecture is the more recent of a series of similar conjectures in different settings. The first conjecture of this series is known as \(p\)-adic monodromy conjecture and it was suggested by Igusa.
Given a polynomial \(f\in\fZ[x_1,\dots, x_n]\), this conjecture suggests a connection between the associated \(p\)-adic zeta function and the local monodromy eigenvalues of \(f\colon \fC^n\to \fC\). After Kontsevich's work on motivic integration,  Denef and Loeser proposed an upgrade of this conjecture in \cite{dl}. Given a polynomial \(f\in \fC[x_1,\dots, x_n]\), they introduced the associated motivic zeta function and formulated the motivic monodromy conjecture for hypersurface singularities, which statement is analogue to the one of the \(p\)-adic monodromy conjecture. 

The main aim of this survey is to present this topic in a rigorous way. In Section \ref{1} we give the definition of a degeneration of a Calabi-Yau varieties and we explain what a model is. In Section \ref{2} we recall some basic definitions in motivic integration and state the monodromy property for Calabi-Yau varieties. The main technical issue is that the motivic zeta function is a power series on a certain localization of an equivariant Grothendieck ring of varieties, so we have to be careful when we talk about poles of the function. For the definition of the monodromy zeta function and the A'Campo formula we refer to \cite{n}. the equivariant version of the Grothendieck ring of varieties was introduced in \cite{h}. The definition of the motivic zeta function and the Denef and Loeser's formula are the ones given in \cite{hn2}. The definitions of poles of a power series with coefficients in the Grothendieck ring of varieties were given in \cite{nx}.
In Section \ref{3} we talk about abelian varieties, Calabi-Yau varieties which admit an equivariant Kulikov model and Kummer surfaces (we refer to \cite{hn2} for the first two cases and to \cite{o} for the last one). In Section \ref{4} we discuss triple-points-free models of K3 surfaces (we refer to \cite{j0}, \cite{j} and \cite{l}). Finally, in Section \ref{5} we briefly discuss  some directions of the research in this topic.

\subsection*{Acknowledgement}
I would like to thank Prof. J. Nicaise for his help during the preparation of this paper.
This work was supported by the Engineering and Physical Sciences Research Council [EP/L015234/1]. The EPSRC Centre for Doctoral Training in Geometry and Number Theory (The London School of Geometry and Number Theory), University College London


\section{Models for Calabi-Yau variety over \(K\)}\label{1}

In the introduction we talked about degenerations of Calabi-Yau varieties, the picture we described was intuitive and geometric, however, it was not rigorous. In this section we define rigorously what we mean by degenerations of Calabi-Yau varieties and what  a model is.

We fix the following notations \(k = \fC\), \(K = \fC((t))\) and \(R = \fC[[t]].\) For every positive integer \(d\), we have \(K(d) = \fC((\sqrt[d]t)), R(d) = \fC[[\sqrt[d]t]].\)
We set \(\overline K = \fC\{\{t\}\},\) the field of Puiseux series:
\[\fC\{\{t\}\} = \bigcup\limits_{d = 1}^\infty R(d);\]
this is an algebraic closure of \(K\).
The group \(\gal( K(d)/ K)\) is canonically isomorphic to \(\mu_d,\) the group of \(d\)-th roots of unity in \(k\). The group \(\widehat\mu\) is the profinite group of roots of unity in \(k\), it is obtained as \(\varprojlim \mu_d\). We have that \(\widehat\mu\cong\gal(\overline K/ K)\).

\begin{definition}[Calabi-Yau varieties]
A Calabi-Yau variety \(X\) over \(K\) is a smooth, proper, geometrically connected variety with trivial canonical sheaf. 
\end{definition}
\begin{definition}[Abelian variety]
An abelian variety \(A\) over \(K\) is a smooth, proper, geometrically connected commutative group-scheme over \(K\). 
\end{definition}
\begin{remark}All abelian varieties are Calabi-Yau varieties. In fact, since they are group varieties, their tangent bundle is trivial. Thus it follows that the top exterior power of the cotangent bundle is the trivial line bundle.
\end{remark}
\begin{definition}[K3 surfaces]
A K3 surface over \(K\) is a 2-dimensional Calabi-Yau variety over \(K\) with \(H^1(X,\sO_X) =0.\) 
\end{definition}
\begin{remark}If we restrict out attention to surfaces over \(\overline K\), it follows  by the Enriques-Kodaira classification, that a 2 dimensional Calabi-Yau variety is either an abelian or a K3 surface.
\end{remark}
\begin{definition}[Model]\label{model}
Let \(X\) be a proper and smooth \(K\)-scheme; a model for \(X\) over \(R\) is a flat \(R\)-algebraic space \(\sX\) endowed with and isomorphism of \(K\)-schemes: \(\sX_K = \sX\times_R K \to X\). We say that \(\sX\) is a strict normal crossing model (snc model) for \(X\) if it is regular and proper over \(R\), and \(\sX_k = \sX\times_R k  \) is a  strict normal crossing divisor on \(\sX\). The special fiber, under this definition, need not to be reduced. If the fiber is reduced then the model is called semistable. The surface \(\sX_k\) is called special fiber, while \(X_K\) is called generic fiber.
\end{definition}
\begin{remark}
In this remark we compare this definition of a model to the one that we gave in the introduction. The space \(\spec R\) corresponds to the disk \(D\); the generic fiber corresponds to the degeneration over the punctured disk, while the special fiber corresponds to the fiber over \(0\). The generic fiber over \(\overline K\) corresponds to the universal fiber of the degeneration, \textit{i.e.} the base change to a universal covering space of \(D^\ast\).
\end{remark}
\begin{remark}\label{sncmodel}
Snc models always exists by Hironaka's resolution of singularities, while semistable models do not, in general. However, by the semi-stable reduction Theorem \cite{kkms}, given any proper model \(\sX\) of \(X\) there exist a positive integer \(d\) and a semistable model \(\sX_d\) of \(X\times_KK(d)\) that dominates \(\sX\times_R R(d)\). If, moreover, \(X\) is projective, then also \(\sX_d\) may be taken projective. It is important to notice that \(\sX_d\) is not a model for \(X\).
\end{remark}

\begin{notation}\label{ei}
Fix an snc model  \(\sX\) for \(X\). The special fiber of this model is \(\sX_k = \sum\limits_{i\in I}N_i E_i\), then for any \(J\subset I\) we define \(E_J = \bigcap\limits_{j\in J} E_j\) and \(E^\circ_J = E_J \setminus \left(\bigcup_{i\not\in J} E_i\right)\).
\end{notation}


\begin{example}\label{quartic}
Let \(X\) be the K3 surface in \(\fP^3_K\) given by the equation:
\[x^2w^2 + y^2 w^2 + z^2 w^2 + x^4 + y^4 + z^4 + tw^4 =0.\]
Let \(\sX\) be the closed subscheme of  \(\fP^3_R\) given by the same equation. The scheme \(\sX\) is regular, but the special fiber \(\sX_k\) is a singular surface. The only singular point is \(P = [0,0,0,1]\); the singularity at that point is a canonical singualirity of type \(A_1\). If we blow up \(\sX\) at \(P\), we obtain an snc model for \(X\), call it \(\sX'\).

The special fiber of this model is non-reduced, so in particular it is not a semistable model. As a divisor the special fiber of this model is \(\sX'_k = D + 2E,\) where \(D\) is the proper transform of \(\sX_k\), which is a smooth K3 surface, and \(E\cong \fP^2_k.\) The intersection of \(D\) and \(E\) is transverse and it is a smooth conic in \(E\).
\end{example}

\section{Monodromy property for Calabi-Yau varieties}\label{2}
In this section we explain what the motivic monodromy property for Calabi-Yau varieties is. First we explain what is the monodromy action; then we introduce the motivic zeta function and explain how they are, conjecturally, related.

\subsection{The monodromy action}
Recall that we fixed  \(\overline K = \fC\{\{t\}\},\) the field of Puiseux series; let \(\sigma\) be the canonical topological generator of the Galois group \(\gal(\overline K/K).\) The generator \(\sigma\) can be described as \(\sigma = \left(\exp\left(\frac{2\pi i}{d}\right)\right)_{d>0}\) and is called the monodromy operator.
\begin{definition}[Monodromy eigenvalue]\label{me}
If \(X\) is a smooth, proper variety over \(K\), then, for all \(i\), the monodromy operator \(\sigma\) acts on the \(l\)-adic cohomology group \(H^i(X\times_K \overline K, \fQ_l)\). We say that \(\lambda\) is a monodromy eigenvalue if there exists an index \(i\) such that \(\lambda\) is an eigenvalue the action of \(\sigma\) on \(H^i(X\times_K \overline K, \fQ_l)\).
\end{definition}
\begin{definition}[Monodromy zeta function]
The monodromy zeta function of \(X\) is defined as
\[\zeta_X(T)= \prod\limits_{i>0}\left(\det\left(T\cdot \id - \sigma|_{ H^i(X\times_K \overline K, \fQ_l)}\right)\right)^{(-1)^{i+1}}\in \fQ_l(T).\]
\end{definition}
\begin{remark}
The monodromy zeta function does not encode all the information about the monodromy eigenvalues, in fact, some cancellations may occur. Moreover, there is an even more natural function that encodes all the monodromy eigenvalue - namely, the product of all the characteristic polynomials. 
One of the main reason why the monodromy zeta function is such a useful tool to study the monodromy eigenvalues is that Theorem \ref{acampo} gives an alternative, and often easier, way to compute this function, while computing all the characteristic polynomials of the monodromy action is usually more complicated.
\end{remark}
\begin{remark}
If \(X\) is a K3 surface over \(K\) then  \(H^i(X\times_K \overline K\, \fQ_l)\) is trivial in odd degrees; as a consequence, the monodromy zeta function of a K3 surface has all the monodromy eigenvalues as poles.
\end{remark}
\begin{theorem}\label{acampo}[A'Campo's formula, \cite{n} Theorem 6.2.6]
Let \(X\) be a smooth, proper \(K\)-variety, and fix an snc model \(\sX\). The special fiber of this model is \(\sX_k = \sum_{i\in I} N_i E_i;\) then the monodromy zeta function is given by
\[\zeta_X(T) = \prod\limits_{i\in I}\left(T^{N_i} - 1\right)^{-\chi_{top}(E_i^\circ)},\]
\(\chi_{top}\) is the topological Euler characteristic and \(E_i^\circ = E_i\setminus\bigcup_{j\neq i} E_j\).
\end{theorem}
\begin{remark}
If follows from A'Campo formula that \(\zeta_X(T)\in \fQ(T)\subset  \fQ_l(T) \).
\end{remark}
\subsection{The motivic zeta function}
Next we discuss the second main ingredient to understand the monodromy property: the motivic zeta function. First we need some background in motivic integration. Since we want to keep track of the Galois action of \(\gal(\overline K/ K)\), instead of working with the usual Grothendieck ring of varieties, we use an equivariant version of it.
\begin{definition}[Category of \(G\)-schemes over \(k\)]
The category \(\sch_{k,G}\) is the category that has:
\begin{enumerate}
\item as objects separated \(k\)-schemes of finite type with good \(G\)-action;
\item as morphisms \(G\)-equivariant morphisms of \(k\)-schemes.
\end{enumerate}
We say that the action of the group \(G\) on the scheme \(X\) is good if \(X\) has a finite partition in \(G\)-stable affine subschemes. 
\end{definition}
\begin{remark}
The definition of good action of the group \(G\) on the scheme \(X\) that  we use is the one given in \cite[(2.2.1)]{hn2}, and it is weaker than the one that is commonly used \textit{i.e.} that \(X\) can be covered by \(G\)-stable affine open subschemes. The reason we preferred our definition is that it can be generalized to algebraic spaces, and it gives rise to the same Grothendieck ring as the usual one.
\end{remark}
\begin{definition}[Equivariant Grothendieck ring, \cite{h} Definition 4.1]
Fix a finite group \(G\). The equivariant Grothendieck ring of \(G\)-varieties over \(k\) is the ring generated as an abelian group by the isomorphism classes of objects \(X\in \sch_{k,G}\), with the ring structure given by the fiber product over \(k\), with the additional relations:
\begin{enumerate}
\item (Scissor relation) Given a \(G\)-scheme \(X\) and a closed \(G\)-subscheme \(Y\), then 
\[[X] = [Y] + [X\setminus Y].\]
\item Given \(A_1\) and \(A_2\) two \(G\)-equivariant affine bundles of rank \(d\) over a scheme \(S\in \sch_{k,G}\) we have
\([A_1] = [A_2].\)
\end{enumerate}We denote this ring by \(K_{0,k}^G.\)
\end{definition}
\begin{remark}
An algebraic space \(X\) with good \(G\)-action defines a class in the equivariant Grothendieck ring. In fact, since \(X\) admits a partition in \(G\)-stable affine schemes, it is possible to use the scissor relation to construct the class in the equivariant Grothendieck ring.
\end{remark}
To our purposes the ring \(K_{0,k}^G\) is not enough, we have to invert an element. Fix \(\fL = [\fA^1]\), where the action of \(G\) is trivial, we define the ring \(\sM^G_k =K_{0,k}^G[\fL^{-1}] \). If we have a profinite group \(\widehat G = \varprojlim G_i,\) with all the groups \(G_i\) finite, then we define:
\[\sM^{\widehat{G}}_k =  \varinjlim\sM^{G_i}_k\]

\begin{definition}[Equivariant weak N\'eron model]\label{weak}
Let \(X\) be a smooth proper \(K\)-scheme. For every \(d>0\), set \(X(d) = X\times_K K(d)\). There is an action of \(\mu_d\) on \(X(d)\). An equivariant weak N\'eron model for \(X(d)\) is a separated and regular algebraic space \(\sX\) over \(R(d)\), with a good action of \(\mu_d.\) Moreover, we require that there exists an isomorphism of \(K(d)\) schemes \(\sX_{K(d)}\to X(d)\) which is \(\mu_d\)-equivariant and such that the natural map \(\sX(R(d))\to X(K(d))\) is a bijection. 
\end{definition}
\begin{remark}An equivariant weak N\'eron model always exists, as explained in \cite[2.2.3]{hn2}. If we have \(\sX\) an snc model for a Calabi-Yau variety \(X\), then it is possible to construct an equivariant weak N\'eron model for \(X(d)\). We have just to normalize the pullback of \(\sX\) on the new base, apply a \(\mu_d\) equivariant resolution of singularities and then restrict to the smooth locus.
\end{remark}
Now we are ready to talk about motivic integration. We fix our variety \(X\) over \(K\) with trivial canonical bundle, and choose a volume form \(\omega\) on \(X\). Denote by \(\omega_d\) the pullback of \(\omega\) to \(X(d)\). Choose a weak equivariant N\'eron model \(\sX\) for \(X(d)\). For every connected component \(C\) of the special fiber \(\sX_k\), the order of \(\omega_d\) along \(C\) is the unique integer \(n\) such that \(t^{-n/d}\omega_d \) is a generator for the sheaf \(\omega_{\sX/R(d)}\) locally at the generic point of \(C\). For every integer \(i\), let \(C(i)\) the union of the connected components of \(\sX_k\) of order \(i\). It is important to remark that \(C(i)\) is stable under the action of \(\mu_d\).
\begin{definition}[Motivic integral]
With the above notation the motivic integral of \(\omega_d\) on \(X(d)\) is defined as:
\[\int\limits_{X(d)}|\omega_d| = \sum\limits_{i\in \fZ}[C(i)]\fL^{-i}\in \sM^{\widehat\mu}_k.\]
\end{definition}
\begin{proposition}[\cite{hn2}, Proposition 2.2.5]
The motivic integral of \(\omega_d\) on \(X(d)\) is independent of the choice of the weak N\'eron model \(\sX.\)
\end{proposition}
\begin{definition}[Motivic zeta function]\label{mzf}
Let \(X\) be a Calabi-Yau variety over \(K\), choose \(\omega\) a volume form on \(X\). The motivic zeta function of \(X\) with respect to the volume form \(\omega\) is:
\[Z_{X,\omega} (T) = \sum\limits_{d>0}\left ( \int\limits_{X(d)}|\omega_d|\right )T^d \in \sM^{\widehat\mu}_k[[T]].\]
\end{definition}
As in the case of the monodromy zeta function, once we fix an snc model, there exists an alternative way to compute the motivic zeta function. Assume that \(\sX\) is an snc model for \(X\), and let \(\sX_k = \sum_{i\in I} N_iE_i\) be the special fiber, the reduced special fiber is \(\sX_{k, red} = \sum_{i\in I} E_i\). The volume form \(\omega\) on \(X\) defines a rational section of the line bundle \(\omega_{\sX/R}(\sX_{k,red} - \sX_k)\) on \(\sX\). To this section we can associate a divisor supported on the special fiber. This divisor is of the form \(\sum_{i\in I} \nu_i E_i.\) By numerical data of \(E_i\) we mean the couple \((N_i, \nu_i)\). For any subset \(J\subset I\) consider the varieties \(E_J\) and \(E_J^\circ\), defined as in Notation \ref{ei}, and set \(N_J = \gcd\{N_j|j\in J\}\). Let \(\sX(N_J)\) be the normalization of \(\sX\times_R R(N_J).\)  Then \(\sX(N_J)\times_\sX E_J^\circ\) is a Galois cover of \(E_J^\circ\) (an explicit description of this cover can be found in Section 2.3 of \cite{n}). Now we are ready to give an alternative description of the motivic zeta function.
\begin{theorem}[Denef and Loeser's formula, \cite{hn2}]
In the above setting we have:
\[Z_{X,\omega}(T) = \sum\limits_{\emptyset\neq J\subset I}[\widetilde{E^\circ_J}](\fL - 1)^{|J| - 1}\prod\limits_{j\in J}\frac{\fL^{-\nu_j}T^{N_j}}{1 - \fL^{-\nu_j}T^{N_j}}\in \sM^{\widehat\mu}_k[[T]]\]
\end{theorem}
It is clear from Definition \ref{mzf} that the motivic zeta function is a power series with coefficients in the ring \(\sM^{\widehat\mu}_k\). However, the definition of poles of this function is not immediate.
\begin{definition}[Poles, \cite{nx}]
Let \(Z(T) \in \sM^{\widehat\mu}_k \left[T, \frac{1}{1 - \fL^a T^b}\right]_{(a,b)\in S}\) be a rational function over \(\sM^{\widehat\mu}_k\), and choose \(q\in \fQ\). We say that the rational number \(q\) is a pole of order at most \(m\geq 1\) for the function \(Z(T)\) if there exists a set \(\sP\) whose elements are multisets contained in \(\fZ\times\fZ_{>0}\) such that:
\begin{enumerate}
\item each multiset \(P\in \sP\) contains at most \(m\) elements \((a,b)\) such that \(\frac{a}{b} = q\) and
\item \(Z(T)\) is an element in the \(\sM^{\widehat\mu}_k [T]\)-submodule of \(\sM^{\widehat\mu}_k [[T]]\) generated by the elements of the form
\[\prod\limits_{(a,b)\in P}\frac{1}{(1 - \fL^a T^b)},\]
for all multisets \(P\in \sP\).\end{enumerate}
The order of a pole at \(q\) is the minimal \(m\) such that \(Z(T)\) has a pole at \(q\) of order at most \(m\).
\end{definition}
\begin{remark}
The reason why the definition of a pole is so involved is that the ring  \(\sM^{\widehat\mu}_k\) is not a domain. 
Actually in \cite{p} it is proven that, if the base field has characteristic zero, not even the usual Grothendieck ring of varieties \(K_{0,k}\) is a domain. 
Moreover, if \(k = \fC\) it was proven in \cite{b} that \(\fL\) is a zero-divisor in \(K_{0,k}\). 
It was proven in\cite{e} that \(\sM_k = K_{0,k}[\fL^{-1}] \) is not a domain, the background to this result is presented in Appendix A of \cite{c}; since \(\sM_k\) injects into \(\sM^{\widehat\mu}_k\) (as varieties with trivial \(\widehat\mu\) action) it follows that the latter is not a domain either.
\end{remark}
At this point we have all the background needed to explain what the monodromy property for Calabi-Yau varieties is.
\begin{definition}[\cite{hn2}, Definition 2.3.5]
Given a Calabi-Yau variety \(X\) over \(K\) and a volume form \(\omega\) on \(X\) we say that the couple $(X,\omega)$ satisfies the monodromy property if there exists a finite set \(S\subset \fZ\times \fZ_{>0}\) such that \(Z_{X,\omega}(T)\) is an element of the sub-ring 
\[\sM^{\widehat\mu}_k \left[T, \frac{1}{1 - \fL^a T^b}\right]_{(a,b)\in S}\subset \sM^{\widehat \mu}_k[[T]],\]
and, for any couple \((a,b)\in S\), we have that \(\exp\left(2\pi \sqrt{-1} \frac{a}{b}\right)\) is a monodromy eigenvalue as in definition \ref{me}; which means that it is an eigenvalue of the action of \(\text{Gal}(\overline K/ K)\) on \(H^i(X\times_K \overline K, \fQ_l)\) for some \(i\geq0\) and every embedding of \(\fQ_l\) in \(\fC\). 
\end{definition}


\section{Some well-understood cases}\label{3}
We are far from a complete understanding of whether or not Calabi-Yau varieties satisfy the monodromy property. It was proven for certain families, under some more restrictive hypothesis on the type of Calabi-Yau variety or on the type of degeneration. 
The conjecture was first proven for abelian surfaces in \cite{hn1}, this proof uses the theory of N\'eron models.
A generalization of this result was obtained in \cite{hn2}; in this paper it was proven that degenerations of Calabi-Yau varieties which admit equivariant Kulikov models satisfy the monodromy property and that abelian varieties admit such models. 
Using the aforementioned result of \cite{hn2}, it was also possible to prove the monodromy property for Kummer surfaces: it was showed in \cite{o} that they admit equivariant Kulikov models. 
All the examples we mentioned have in common that the motivic zeta function has a unique pole. This is not by chance, in fact, we have the following, more general, result.

\begin{theorem}[\cite{hn2}, Theorem 3.3.3]\label{onepole}
Let \(X\) be a Calabi-Yau variety of dimension \(n\) over \(K\), with a volume form \(\omega\). Choose an snc model for \(X\), let \(\sX_k = \sum_i N_i E_i\) be the special fiber of this model and denote with  \(\nu_i\) the vanishing order of \(\omega\) on \(E_i\). Let \(\min(\omega) = \min_i \frac{\nu_i}{N_i}\).
Under these hypothesis one of the eigenvalue of the action of \(\sigma\) on \(H^{n}(X\times_K \overline K, \fC)\) is \(\exp(-2\pi i\min(\omega))\).
\end{theorem}
Indeed, in all the cases we mentioned at the beginning of the section, it was proven that the motivic zeta function of \((X,\omega)\) has a unique pole at \(1 - \min(\omega)\), and since \(\exp(-2\pi i\min(\omega))\) is a monodromy eigenvalue, then the monodromy property holds.
In Section \ref{4} we present a family of degenerations of K3 surfaces that satisfies the monodromy property, but whose motivic zeta functions may have more than one pole, most of these results are the work of \cite{j0} and \cite{j}.

\begin{theorem}[\cite{hn2}, Theorem 4.2.2]\label{monab}
Fix an abelian variety \(A\) with a volume form \(\omega\). Choose an snc model \(\sA\), and let \(\sA_k = \sum_i N_i E_i\) be the special fiber. We denote with \(\nu_i\) the vanishing order of \(\omega\) on \(E_i\). The motivic zeta function \(Z_{A,\omega}(T)\) has a unique pole at \(q = 1 - \min(\omega)\). More precisely:
\[Z_{A,\omega}(T)\in\sM^{\widehat\mu}_k \left[T, \frac{1}{1 - \fL^a T^b}\right]_{(a,b)\in S;\; a/b = q}\subset \sM^{\widehat \mu}_k[[T]].\]
\end{theorem}
It follows that the monodromy property holds for abelian varieties from Theorem \ref{monab} and Theorem \ref{onepole}.

Now we are ready to introduce the definition of equivariant Kulikov models for Calabi-Yau varieties.
\begin{definition}[Equivariant Kulikov model]
Let \(X\) be a Calabi-Yau variety over \(K\), and fix a positive integer \(d\). A Kulikov model for \(X\) over \(R(d)\) is a regular, proper and flat algebraic space \(\sX\) over \(R(d)\) such that:
\begin{enumerate}
\item there is an isomorphism of \(K(d)\)-schemes:
\[\xymatrix{\sX_{K(d)}\ar[r]^{\sim\quad\;\;}\ar[dr]&X\times_KK(d)\ar[d]\ar[r]&X\ar[d]\\&\spec(K(d))\ar[r]&\spec K},\]
\item the special fiber \(\sX_k\) is a divisor with normal crossing
\item the logarithmic relative canonical bundle \(\omega_{\sX/R(d)}(\sX_{k, red} - \sX_k)\) is trivial.
\end{enumerate}
 We say that the Kulikov model \(\sX\) is equivariant if the Galois action of \(\mu_d\) on \(X\times_K K(d)\) extends to \(\sX\).
\end{definition}
\begin{remark}
Not all snc models are Kulikov models; for instance, the model we constructed in Example \ref{quartic} is not a Kulikov model over \(K\). It is easy to check that the K3 surface of Example \ref{quartic} admits a Kulikov model over \(K(2)\), however, we will show in Example \ref{quartic2} that this model is not equivariant. Indeed, in Example \ref{quartic2} we show something stronger, \textit{i.e.} that the K3 surface of Example \ref{quartic} does not admit an equivariant Kulikov model for any \(d\).
\end{remark}
\begin{theorem}[\cite{hn2}, Theorem 5.3.2]\label{monekm}
Let \(X\) be a Calabi-Yau variety, with a volume form \(\omega\). Assume that \(X\) admits an equivariant Kulikov model over \(R(d)\) for some positive \(d\). Then the motivic zeta function of \((X, \omega)\) has a unique pole at \(q = 1 - \min_i \frac{\nu_i}{N_i} \).
\[Z_{X,\omega}(T)\in\sM^{\widehat\mu}_k \left[T, \frac{1}{1 - \fL^a T^b}\right]_{(a,b)\in S, a/b = q}\subset \sM^{\widehat \mu}_k[[T]].\]
\end{theorem}
\begin{remark}
As in the previous case, the monodromy property for Calabi-Yau varieties admitting an equivariant Kulikov model follows from Theorem \ref{monekm} and Theorem \ref{onepole}. This result generalize Theorem \ref{monab}, in fact it was proven that abelian varieties admit an equivariant Kulikov model in \cite[Theorem 5.1.6]{hn2}.
\end{remark}
An other consequence of Theorem \ref{monekm} is that the monodromy property holds for Kummer surfaces, this was proven in \cite{o}. 
\begin{definition}[Kummer surface]
Let \(A\) be an abelian surface over \(K\), consider an involution \(\iota\), and call \(A_\iota\) the fixed point of \(\iota\). Let \(\widetilde A\) be the blow of \(A\) at \(A_\iota\); the involution \(\iota\) acts regularly on  \(\widetilde A\). Call \(X\) the quotient of \(\widetilde A\) by the action of \(\iota\); \(X\) is a smooth K3 surface over \(K\). Any K3 surface that can be obtained in this way is called Kummer surfaces.
\end{definition}
\begin{theorem}[\cite{o}, Theorem 6.2]\label{monk}
Let \(X\) be a Kummer surface, then there exists a minimal \(d_0 >0\) such that \(X(d_0)\) has an equivariant Kulikov model, moreover, if \(d>0\) is such that  \(X(d)\) admits an equivariant Kulikov model, then \(d_0|d\). In particular, Kummer surfaces satisfy the monodromy property.
\end{theorem}
Theorem \ref{onepole} was also used in \cite{hn2} to prove the monodromy property for some Calabi-Yau varieties which were not in any of these families. As far as we know, the cases we listed in this section are the only ones where Theorem \ref{onepole} was used to prove the monodromy property.


\section{Triple-Points-Free models of K3 and why they are interesting}\label{4}
In this section we discuss triple-points-free models of K3 surfaces. If a K3 surface admits a triple-points-free model, it satisfies the monodromy property, however, the motivic zeta function of this K3 surface may have more then a single pole. In particular K3 surfaces admitting a triple-points-free model may not admit an equivariant Kulikov model. 

In \cite{cm} Crauder and Morrison described special fiber of relatively minimal, triple-points-free snc models of surfaces with trivial pluricanonical bundles. Many additional results to describe the combinatorics of the special fiber were given by Jaspers in \cite{j}.
The motivic zeta function and the monodromy property for these surfaces were studied in \cite{j0} and \cite{j}. In Corollary 4.2.4 of \cite{j} it was shown that poles of the motivic zeta function may be recovered from the combinatorial data of the model, in particular it was shown that besides the minimal pole, there are additional poles as soon as the triple-points-free model has a so-called conic flower. 

Regarding the monodromy property, it was proven that under some additional conditions on the special fiber, K3 surfaces admitting a triple-points-free model satisfy the monodromy property. 
In  Appendix B of \cite{j} there is a classification of the possible special fibers of triple-points-free models that do not satisfy the monodromy property; such surfaces are known as combinatorial countercandidates. In \cite{l} we proved that these combinatorial countercandidates do not exist. As a consequence, the monodromy property holds for K3 admitting triple-points-free models. 
\begin{definition}[Triple-points-free model]
Given a K3 surface \(X\), a triple-points-free model \(\sX\) of \(X\) is a relatively minimal snc model such that given any three distinct irreducible components \(E_i, E_j, E_k\) of the special fiber \(\sX_k\), 
\[E_i\cap E_j\cap E_k = \emptyset.\]
Since there are not triple intersections, the dual complex of the special fiber of this model is a graph, we call it \(\Gamma\); to each vertex of \(\Gamma\) we associate the weight 
\[\rho_i = \frac{\nu_i}{N_i} + 1;\]
we denote by \(\Gamma_{min}\)  the subgraph of \(\Gamma\) of components of minimal weight. 
\end{definition}

There is a very explicit description of possible special fibers of triple-points-free models for K3 surfaces. 
\begin{theorem}[Crauder-Morrison classification for K3 surfaces, \cite{cm} and \cite{j}] 
Let X be a smooth, proper \(K3\) surface over \(K\), let \(\sX\) be a relatively minimal triple-points-free model of \(X\), then \(\sX\) has the following properties:
\begin{enumerate}
\item \(\Gamma_{min}\) is a connected subgraph of \(\Gamma\). It is either a vertex or a chain. 
\item Each connected component of \(\Gamma\setminus \Gamma_{min}\) is a chain (called flower) \(F_0-F_1-\cdots-F_l\) where only \(F_l\) meets \(\Gamma_{min}\). The weights strictly decrease along these flowers, \(F_0\) being the one with maximal weight. The surface \(F_0\) is either minimal ruled or isomorphic to \(\fP^2\). If it is isomorphic to \(\fP^2\) then \(F_0\cap F_1\) is either a line or a conic. The other components are minimal ruled surfaces, and \(F_i\cap F_{i + 1}\) and \(F_i \cap F_{i - 1}\) are both sections of the ruling.
\item If \(\Gamma_{min}\) is a single vertex,  there are three possible cases for the corresponding surface. It is either isomorphic to \(\fP^2\) or it is a ruled surface or the canonical divisor is numerically trivial. If \(\Gamma_{min}\) is a point, we call the model a flowerpot degeneration.
\item If \(\Gamma_{min}\) is a chain \(V_0 - V_1 - \cdots - V_{k+1}\), then we can describe the components of the chain. If \(i\neq 0, k + 1\), \(V_i\) is an elliptic ruled surface, and \(V_{i - 1}\cap V_i\) and \(V_i \cap V_{i + 1}\) are both sections of the ruling; if \(i = 0, k + 1\), then \(V_i\) is isomorphic to \(\fP^2\), or  it is a, rational or elliptic, ruled surface. If \(\Gamma_{min}\) is a chain we call the model a chain degeneration.
\end{enumerate}
\end{theorem}
\begin{remark}
In \cite{cm} there is an even stronger result. Indeed, they classified the special fiber of triple-points-free models with generic fiber \(\sX_K\) with trivial pluricanonical bundle. In \cite{cm} there is also a complete classification of the possible flowers, divided in 21 combinatorial classes.  
If \(\sX_K\) is a K3 surface, this classification was refined in Chapter 3 of \cite{j}.
\end{remark}
One of the main results of \cite{j} is the description of the motivic zeta functions of these models. This description is extremely explicit and using it, it was possible to study the poles. It turned out that poles of the motivic zeta function are closely related to the presence of conic flowers (\textit{i.e.} \(F_0\) is isomorphic to \(\fP^2\) and \(F_0\cap F_1\) is a conic).
\begin{theorem}[\cite{j}, Theorem 4.3.8]\label{amt}
Let \(X\) be a a K3 surface over \(K\) with a volume form \(\omega\). Choose a triple-points-free model \(\sX\), whose special fiber is \(\sX_k = \sum_i N_i E_i,\) with numerical data \((N_i,\nu_i)\).
Then \(q\in \fQ\) is a pole of \(Z_{X,\omega}(T)\) if and only if there exists an \(i\) such that the numerical data of \(E_i\) satisfy \(q = -\nu_i/N_i\) and such that:
\begin{enumerate}
\item either \(\rho_i\) is minimal
\item or \(E_i\) is the top of a conic flower.
\end{enumerate} 
Moreover, while in the second case the pole is always of order 1, in the first one it is of order 1 if \(\sX\) is a flowerpot degeneration and of order 2 if it is a chain degeneration.
\end{theorem}
\begin{remark}
In \cite[Appendix A]{j} there is a Python code that describes the contribution of each flower to the motivic zeta function. 
\end{remark}
\begin{remark}
In \cite{j} it was also discussed whether or not a K3 surface \(X\) admitting a triple-points-free model \(\sX\) satisfies the monodromy property. In \cite[Theorem 5.2.1]{j} it was proven that flowerpot degenerations satisfy the monodromy property. In \cite[Theorem 5.3.1]{j} it was proven that chain degenerations with some extra assumptions satisfy the monodromy property. However, it was not clear if those assumptions were enough to prove the monodromy property for any chain degeneration. 

The strategy adopted in \cite{j} was to try to prove the conjecture by contradiction. 
The assumption was that there exists a chain degeneration that does not satisfy the monodromy property, this mean that some poles of the motivic zeta function of this degeneration do not correspond to monodromy eigenvalues. 
The first step was to use A'Campo's formula and Denef and Loeser's formula to deduce some information on the geometry of the special fiber of this triple-points-free model. The last step would have been to prove that these surfaces - called combinatorial countercandidates - do not appear as the special fiber of any triple-points-free model. This would have proven the monodromy property for all K3 surfaces admitting a triple-points-free model.

In \cite[Section 6]{j} there is a description of some topological properties of these combinatorial countercandidates, and this description become even more explicit in \cite[Appendix B]{j}, where all the possible countercandidates are listed.
The kind of information we have about the numerical countercandidates is:
\begin{enumerate} 
\item the central fiber of the degeneration is a chain degeneration, the chain is \[V_0 - V_1 - \cdots - V_{k+1}.\] 
\item The surface \(V_0\) satisfies the following properties: 
\begin{enumerate}
\item it contains a smooth elliptic curve \(D\) such that \(D = -K_{V_0}\);
\item it is a smooth rational ruled non-minimal surface obtained from an Hirzebruch surface by \(l_0\) blow-ups whose centers lies in the image of \(D\);
\item it contains at least \(h\) smooth rational curves of self intersection -2, we denote them by \(C_i\). The curves \(C_i\) are all disjoint, moreover from the adjunction formula it follows that they do not intersect \(D\);
\end{enumerate}
\item the other internal components of the chain are elliptic or rational ruled surfaces, and they are obtained from minimal surface by some blow-up along the intersection of two components. Moreover, we know that there is a certain number of -2 rational curves, disjoint from the sections \(V_i\cap V_{i+1}\).
\end{enumerate}
Whether or not any combinatorial countercandidate exists was left as an open problem in \cite{j}. 
It was proven in \cite{l}  that they do not exist, all the information required to prove this result are the one regarding the surface \(V_0\).
\end{remark}
\begin{theorem}[\cite{j}, \cite{l}]
Let \(X\) be a K3 surface that admits a triple-points-free model \(\sX\). Then the monodromy property holds for \(X\).
\end{theorem}
To conclude this section we come back to the quartic surface of Example \ref{quartic}. 
\begin{example}[\cite{hn2}, \cite{j0},\cite{j}]\label{quartic2}
Let X be the K3 surfaces in \(\fP^3_K,\) given by the equation:
\[x^2w^2 + y^2 w^2 + z^2 w^2 + x^4 + y^4 + z^4 + tw^4 =0.\]
From Example \ref{quartic} we already know that there exists a model \(\sX'\) such that the special fiber is \(\sX'_\fC = D + 2E,\) where \(D\) is a K3 surface and \(E\cong \fP^2_\fC.\) The curve \(C = D\cap E\) is a smooth conic in \(E\); thus \(E\) is a conic flower of type 2B.
We can now choose a volume form on X (for instance the natural volume form induced by the embedding in \(\fP^3\)). This form extends to a relative volume form on \(\sX'\). We can now compute the numerical data of this relative volume form and we obtain \(\nu_D = 0\) and \(\nu_E = 1.\) The motivic zeta function of this K3 is:
\[Z_{X,\omega} (T) = [\widetilde D^\circ] \frac{T}{1-T} + [\widetilde{E}^\circ] \frac{\fL^{-1}T^2}{1 - \fL^{-1}T^2} + [\widetilde C]\frac{\fL^{-1}T^3}{(1-T)(1 - \fL^{-1}T^2)}. \]
The motivic zeta function of this K3 has two poles, namely \(0\) and \(\frac 1 2\).
\end{example}
\begin{remark}
The quartic surface of Example \ref{quartic2} does not admit any equivariant Kulikov model. This follows from Theorem \ref{monekm}, since its motivic zeta function has two poles.
\end{remark}


\section{What's next}\label{5}
In this section we describe a possible approach to future research in this area. Our idea is to construct further examples of K3 surfaces whose motivic zeta function has more than one pole. Once we have done this, the natural question is whether or not these K3 surfaces satisfy the monodromy property.

From the results of Theorem \ref{amt} it is clear that, in the case of triple-points-free models, additional poles come only from conic flowers. 
However, it is not immediately clear what distinguishes these flowers geometrically. 
To obtain some understanding of this, we tried to study what happens to the model if we contract some flowers. The idea of contracting flowers already appears in \cite{cm1}, in this paper there is also a list of the singularities caused by the contraction of some family of flowers.

From the computations of Example \ref{quartic2}, we can see that the contraction of the flower of type 2B gives a regular model with a special fiber which is irreducible but with a singular point. The situation with non-conic flowers is different, L. Halle showed me that, up to a finite base change, it is possible to contract them smoothly - by this we mean that the resulting model is an snc model. This suggested to us that there might be a relation between poles of the motivic zeta function and singularities of the models.

Assume now that \(X\) is a K3 surface over \(K\), \(\omega\) a volume form on \(X\) and let \(\sX\) be a model for \(X\). To \(X\) we can associate a motivic zeta function \(Z_{X,\omega}\); however,  it may not be possible to use Denef and Loeser's formula with the model \(\sX\). In general it is not regular, and even if it were, the irreducible components of the special fiber \(\sX_k\) might be singular. Of course, by Hironaka's resolution of singularities it is possible to construct an snc model \(\sY\) which dominates \(\sX\), but what we would  like to understand is whether or not it is possible to deduce the presence of additional poles of \(Z_{X,\omega}\) from the singularities of \(\sX\). 
This problem, however, is very generic, and thus pretty tough to approach. So we had to restrict the class of singularities we would consider.

\begin{definition}[Rational double points]
Given \(X\) a normal surface over \(\fC\), we say that a point \(x\in X\) is a rational double point (or ADE or Du Val singularity) if it is a canonical singularity. Given a rational double point singularity, we have that \'etale locally around the point the surface is isomorphic to the closed subset of \(\fA^3\) given by one of the following equations:
\begin{enumerate}
\item \(x^2 + y^2 + z^{n + 1}=0\); these are \(A_n\) singularities;
\item \(x^2 + y^2z + z^{n - 1}=0\), with \(n>3\); these are \(D_n\) singularities;
\item \(x^2 + y^3 + z^{4}=0\); this is the \(E_6\) singularity;
\item \(x^2 + y^3 + z^{3}y=0\); this is the \(E_7\) singularity;
\item \(x^2 + y^3 + z^{5}=0\); this is the \(E_8\) singularity.
\end{enumerate}
\end{definition}

Fix a K3 surface \(X\) which admits a  regular model \(\sX\) whose special fiber \(\sX_k\) has some rational double points (for instance, the K3 surface of Example \ref{quartic2}). Assume, furthermore, that none of the rational double points lie the intersection of some irreducible components of the special fiber and call the set of rational double points \(S\).
Then we can construct an snc model \(\sY\) for \(X\), with a morphism \(\pi\colon \sY\to \sX\)  blowing-up \(\sX\) at some smooth points and at some smooth curves contained in \(\sX_k\). Once we have constructed such an snc model, it is possible to compute the motivic zeta function using Denef and Loeser's formula. 
We can split the motivic zeta function \(Z_{X,\omega}\) in two parts, the one which depends on the strata contained in \(\pi^{-1}(\sX_k \setminus S) \)  and the one which depends on the strata contained in \(\pi^{-1}(S)\). We denote the second part by \(Z_{X,\omega,sing}\), and we refer to it as the contribution of the singularities to \(Z_{X,\omega}\).


\begin{example}
For instance, in Example \ref{quartic2}, the contribution of the \(A_1\) singularity is: 
\[Z_{X,\omega,sing}=[\widetilde{E}^\circ] \frac{\fL^{-1}T^2}{1 - \fL^{-1}T^2} + [\widetilde C]\frac{\fL^{-1}T^3}{(1-T)(1 - \fL^{-1}T^2)};\]
since the motivic zeta function in this example has more than one pole, then \(X\) can't have good reduction over \(K\); in particular we have that the action of \(\gal(\overline K/ K)\) is not trivial.
\end{example}

We are studying the contributions of ADE singularities to the motivic zeta function. In the last part of this paper we briefly discuss some of the questions we would like to answer. Question \ref{q1}, \ref{q2}, and \ref{q3} may be consider short term goals, while Question \ref{q4}, is much wilder, and should be consider a long term goal.
\begin{question}\label{q1}
Assume that \(X\) is a K3 surface that admits a regular model whose special fiber has rational double points, is it true that the motivic zeta function of \(X\) has at least two poles? If not, for which singularities is this true?

A positive answer to this question would 
provide more examples of Calabi-Yau varieties whose motivic zeta function has multiple poles.\end{question}
\begin{question}\label{q2}
In the above setting, can we describe explicitly the contributions of the various singularities?
\end{question}
Having an explicit description of the motivic zeta function  would be not only interesting by itself, but it would be also help understanding something more about the monodromy property for K3 surfaces. This brings us to the next question.
\begin{question}\label{q3}
In the above setting,  does \(X\) satisfy the monodromy property? If not, which obstructions are encountered? 
\end{question}
Requiring that the model \(\sX\) is regular is a strong assumption; we would like to relax this hypothesis. The natural class of threefold singularities that extends Du Val singularities is known as compound Du Val singularities.
\begin{definition}[Compound Du Val singularities]
Given \(X\) a normal threefold over \(\fC\), we say that a point \(x\in X\) is a compound Du Val singularity if for some general hyperplane section \(H\) through \(x\) we have that \(x\in H\) is a Du Val singularity.
\end{definition}
\begin{question}\label{q4} Assume that \(X\) is a K3 surface that admits a model \(\sX\)  whose singularities are compound Du Val singularities. What can we say of the analogs of Questions \ref{q1}, \ref{q2}, \ref{q3} under these least restrictive hypothesis?
\end{question}

A complete answer to \ref{q4} might lead to an improvement of the following result from \cite{lm}. 
\begin{theorem}[\cite{lm}, Theorem 6.1]\label{lm}
Let \(X\) be a K3 surface over \(K\), assume that the action of \(\gal(\overline K/ K)\) on \(\het\) is trivial, then:
\begin{enumerate}
\item there exists a projective model of \(X\) whose special fiber is an irreducible surface with at worst rational double points, whose minimal resolution is a K3 surface;
\item \(X\) has potential good reduction.
\end{enumerate}
\end{theorem}






\end{document}